\documentclass[11pt]{amsart}

\usepackage[T1]{fontenc}
\usepackage{lmodern}
\usepackage[margin=0.9in]{geometry}
\usepackage{amsmath,amssymb,amsthm,mathtools}
\usepackage{booktabs,array,microtype}
\usepackage{needspace}
\usepackage{longtable}
\usepackage[colorlinks=true,linkcolor=blue,citecolor=blue,urlcolor=blue]{hyperref}
\hypersetup{pdftitle={Tameness of the Anick automorphism and the Freudenburg\textendash Moser-Jauslin involution}}
\newtheorem{theorem}{Theorem}
\newtheorem{lemma}{Lemma}
\newtheorem{corollary}[theorem]{Corollary}
\newcommand{\Aut}{\operatorname{Aut}}
\newcommand{\TAut}{\operatorname{TAut}}

\newcommand{\ZZ}{\mathbb Z}
\allowdisplaybreaks[1]
\title[Tameness of Anick and Freudenburg--Moser-Jauslin maps]{Tameness of the Anick automorphism\\ and the Freudenburg--Moser-Jauslin involution}

\author{Alexei Miasnikov$\dag$ and Ivan Shestakov$\ddag$}
\address{Schaefer School of Engineering \& Science, Department of Mathematical Sciences, Stevens Institute of Technology, Castle Point on Hudson, Hoboken, 
NJ 07030-5991, USA}
\email{ amiasnikov@gmail.com} 
\address{ Shenzhen International Center for Mathematics, Shenzhen, China\\ 
Instituto de Matem\'{a}tica e Estat\'{i}stica, Universidade de S\~{a}o
Paulo, Brazil }
\email{ ivan.shestakov@gmail.com}

\thanks{${}^{\dag}$ Supported by the Nebius Academy research grant: AI in Math and Math Education}
\thanks{${}^{\ddag}$ Supported by  IMC SUSTech, Shenzhen, China, and by the grants CNPq  305196/2024-3  and FAPESP 2024/14914-9, Brazil}
\date{September 16, 2026}

\begin{document}

\vspace{1.7em}
\begin{abstract}

Let $\delta=xz-yt$. We give an explicit tame factorization of the Anick automorphism
$$
f=(x+\delta t,y+\delta z,z,t)
$$ 
of $F[x,y,z,t]$. After a linear change of
coordinates, the factorization consists of sixteen elementary
automorphisms, all fixing $z$.
A complete list of the
elementary factors is included in Section~\ref{sec:factors}. The explicit
polynomial factorization identities are defined over $\ZZ[1/2]$; no conclusion in characteristic $2$ is asserted.

We also prove that the Freudenburg--Moser-Jauslin involution
$$
 \tau(a,b,x,y)=
 \bigl(b,a,-b^3x+(1+ab+a^2b^2)y,(1-ab)x+a^3y\bigr).
$$
of $F[a,b,x,y]$ is tame over every field of characteristic different from
$2$ and $3$, and list an explicit factorization in
Appendix~\ref{app:psi}.
\end{abstract}

\maketitle
\section{Statement and conventions}

Let $F$ be a field of characteristic different from $2$, and $F[X_n]=F[x_1,\ldots,x_n]$ the polynomial ring in $n$ variables over $F$.
In the present paper, we discuss the structure of the automorphism group $\Aut(F[X_n])$ of $F[X_n]$ over $F$ for $n\leq 4$.
Let $\varphi:F[X_n]\rightarrow F[X_n]$ be an endomorphism of $F[X_n]$ over $F$. We identify $\varphi$ with the $n$-tuple $(f_1,\ldots,f_n)$ of elements of $F[X_n]$, where $f_i=\varphi(x_i)$ for each $i$.
Then, $\varphi$ is an automorphism if and only if the $F$-algebra $F[X_n]$ is generated by $f_1,\ldots,f_n$.

An automorphism is called {\em elementary} if it has the form
$$
E(i,c,p) =(x_1,\ldots,x_{i-1},c x_i+p,x_{i+1},\ldots,x_n).
$$
where $0\neq c\in F$ and $p$ a polynomial in $F[X_n]$ not containing $x_i$. The
inverse of $E(i,c,p)$ is an elementary automorphism $E(i,c^{-1},-c^{-1}p)$.
For $c=1$ we will write $E(i,1,p)=E_i(p)$ or $E_{x_i}(p)$.
The subgroup of $\Aut(F[X_n])$ generated by all elementary automorphisms is called {\em the tame automorphism group} and is denoted as $\TAut(F[X_n])$.
 \smallskip

It is a fundamental question in polynomial ring theory whether
$\TAut(F[X_n])=\Aut(F[X_n])$ holds for each $n$. 
 The equality is obvious if $n=1$. This also holds
true if $n=2$. It was shown by Jung \cite{Jung} in 1942 when
$F$ is of characteristic zero, and by van der Kulk \cite{vdKulk} in 1953 when
$F$ is an arbitrary field. 

\smallskip
In 1972 M.~Nagata \cite{Nagata} suggested that the following automorphism is not tame in $F[x,y,z]$:
$$
\sigma_{\alpha}(x,y,z)=(x+\alpha z,y+2\alpha x+\alpha^2z,z),
$$
where $\alpha=x^2-yz$. It is easy to see that $\sigma_{\alpha}(\alpha)=\alpha$ and $\sigma_{\alpha}^{-1}=\sigma_{-\alpha}$.

\smallskip
In characteristic zero, the conjecture of Nagata was confirmed by I.~Shestakov and U.~Umirbaev in \cite{SU1,SU2}. The tame-generators question remains open for
$n\geq 4$, even in characteristic zero. In characteristic zero, M.~Smith showed \cite{Smith} that the extension of the Nagata automorphism to the four generators $x,y,z,t$ by setting $\sigma_{\alpha}(t)=t$ becomes tame.
D.Anick (see \cite{Cohn}) suggested the following candidate to be a non-tame automorphism for $n=4:$
$$
 f=(x+\delta t,\ y+\delta z,\ z,\ t),\qquad \delta=xz-yt.
$$

\Needspace{8\baselineskip}
\begin{theorem}\label{thm:main}
Let $F$ be a field of characteristic different from $2$.
The Anick automorphism $f$ is tame. After a linear change of
coordinates, it is a product of sixteen elementary automorphisms,
all fixing $z$, as listed in Section~\ref{sec:factors}.
The explicit polynomial factorization identities are defined over
$\ZZ[1/2]$. In particular, $f$ is tame over $F[z]$.
\end{theorem}

Its tame factorization is explained by three constructions: a determinant
invariant, a splitting of left and right unipotent matrix actions,
and an auxiliary-variable commutator. The sixteen elementary factors
are the expansion of two short commutator words.

\medskip

{\bf Theorem~\ref{thm:fmj}} settles the tameness question raised by Lamy and
Przytycki \cite[p.~608]{LP} for the Freudenburg--Moser-Jauslin involution
$\tau$:
$$
 \tau(a,b,x,y)=
 \bigl(b,a,-b^3x+(1+ab+a^2b^2)y,(1-ab)x+a^3y\bigr)
$$
of $F[a,b,x,y]$ which they suggested might be non-tame. We prove that $\tau$ and
the associated automorphism $\psi$ are tame over every field $F$ of
characteristic different from $2$ and $3$. An explicit factorization of
$\psi$, with $76$ elementary factors taken with respect to three linear
coordinate systems and all fixing $b$, is listed in
Appendix~\ref{app:psi}.

Lamy and Przytycki proved that finite subgroups of the three-variable
tame automorphism group are linearizable in characteristic zero
\cite[Corollary~1.2]{LP}. The nonlinearizable Freudenburg--Moser-Jauslin
$S_3$-action \cite{FMJ}, proved tame in Sections~\ref{sec:fmj}
and~\ref{sec:S3}, shows that this property fails in dimension four.
Thus a finite group of tame automorphisms need not admit a single
polynomial change of coordinates making all its elements linear.

\medskip
The general question remains open: is every automorphism
of $F[X_n]$ tame when $n\geq 4$?
Notice that this is the case for free metabelian groups \cite{BM,R} and this question is open for free metabelian Lie algebras \cite{UU}.

In characteristic zero, a related question is whether $\Aut(F[X_n])$
is generated by linear automorphisms together with automorphisms of
the form $\exp(D)$, where $D$ is a locally nilpotent derivation of
$F[X_n]$. Every elementary  automorphism  $E_i(p)$ has this form:
\[
 E_i(p)=\exp\!\left(p\frac{\partial}{\partial x_i}\right),
 \qquad p\in F[x_1,\ldots,\widehat{x_i},\ldots,x_n].
\]
An elementary automorphism $E(i,c,p)$ is a product of such a shear
and a linear automorphism. Linear automorphisms cannot simply be
omitted, since every locally nilpotent exponential has Jacobian
determinant $1$.

An automorphism is called \emph{stably tame} if its extension fixing
finitely many additional variables is tame. In characteristic zero, the Lie algebra construction
in Section~4 (see also \cite{Smith}) gives a sufficient condition for stable tameness: if $D$
is a triangular locally nilpotent derivation and $a\in\ker D$, then
$\exp(aD)$ becomes tame after adjoining one variable. For the Anick
automorphism, the splitting allows us to use variables already present
in the algebra. These constructions motivate the broader question
whether every polynomial automorphism is stably tame.
\medskip

We use the convention
\[
 (AB)(p)=A(B(p)),\qquad [A,B]=ABA^{-1}B^{-1}.
\]
Thus the rightmost factor acts first on a polynomial. If
$A_{\mathrm{pt}}$ and $B_{\mathrm{pt}}$ denote the corresponding
point maps, then
\[
 (AB)_{\mathrm{pt}}=B_{\mathrm{pt}}\circ A_{\mathrm{pt}}.
\]
Consequently, to evaluate the point map represented by the product
in the table, one performs the coordinate updates in the listed order,
from row $1$ to row $16$.

\section{A determinant-dependent row operation}
Let us write our automorphism in a  matrix form.
Set
\[
 M=\begin{pmatrix}x&y\\t&z\end{pmatrix},\qquad
 U(s)=\begin{pmatrix}1&s\\0&1\end{pmatrix}.
\]
Then $\det M=\delta$, and
\[
 f(M)=U(\delta)M.
\]
Thus $f$ adds the determinant times the second row to the first.
Since the determinant is invariant, the opposite row operation gives
the inverse. The coefficient $\delta$ depends on the changed
generators, so the row-operation description does not by itself
give an elementary polynomial automorphism.

Put $h=\delta/2$ and define
\[
 P(M)=U(h)MU(h)^{\mathsf T},\qquad
 Q(M)=U(h)MU(-h)^{\mathsf T}.
\]
Both maps preserve the determinant, and hence $h$. Under our
composition convention,
\begin{align*}
 PQ(M)
 &=P\bigl(U(h)MU(-h)^{\mathsf T}\bigr)\\
 &=U(h)\,P(M)\,U(-h)^{\mathsf T}=U(2h)M=f(M).
\end{align*}
The right actions cancel; the left parameters add.

Introduce the symmetric and antisymmetric off-diagonal coordinates
\[
 v=\frac{y+t}{2},\qquad u=\frac{y-t}{2},\qquad
 y=v+u,\quad t=v-u.
\]
In the order $(x,v,z,u)$, define
\[
 q_+=xz-v^2,\quad q_-=xz+u^2,\quad
 \Delta=q_++u^2=q_--v^2=\delta.
\]
To write $P$ and $Q$ in these coordinates, let $g$ be any polynomial
and multiply out:
\[
 U(g)MU(g)^{\mathsf T}=\begin{pmatrix}x+g(y+t)+g^2z&y+gz\\ t+gz&z\end{pmatrix},
 \qquad
 U(g)MU(-g)^{\mathsf T}=\begin{pmatrix}x-g(y-t)-g^2z&y+gz\\ t-gz&z\end{pmatrix}.
\]
In the first product $y$ and $t$ both increase by $gz$; hence $v$
increases by $gz$, $u$ is unchanged, and, as $y+t=2v$, the corner entry
is $x+2gv+g^2z$. In the second product $y$ increases and $t$ decreases
by $gz$; hence $u$ increases by $gz$, $v$ is unchanged, and, as
$y-t=2u$, the corner entry is $x-2gu-g^2z$. For polynomials
$h,k\in F[x,v,z,u]$ we therefore define the endomorphisms
\begin{align*}
 T_h&=(x+2hv+h^2z,\ v+hz,\ z,\ u),\\
 R_k&=(x-2ku-k^2z,\ v,\ z,\ u+kz).
\end{align*}
Thus $T_h$ replaces the entries of $M$ by the corresponding entries of
$U(h)MU(h)^{\mathsf T}$, and $R_k$ replaces them by those of
$U(k)MU(-k)^{\mathsf T}$; the parameters $h$ and $k$ are fixed
polynomials and are not themselves transformed. With $h=k=\delta/2$,
written $\Delta/2$ in the new coordinates, this gives
\[
 P=T_{\Delta/2},\qquad Q=R_{\Delta/2}.
\]
The first flow fixes $u$; the second fixes $v$. This is the reason
the splitting is useful: each part has an available auxiliary
coordinate, although the whole automorphism moves both $u$ and $v$.

\begin{lemma}\label{lem:flows}
\begin{enumerate}
\item For every polynomial $h$, the endomorphism $T_h$ fixes $q_+$, $u$
and $z$. For every polynomial $k$, the endomorphism $R_k$ fixes $q_-$,
$v$ and $z$.
\item If $h'\in F[q_+,u,z]$, then $T_hT_{h'}=T_{h+h'}$ for every
polynomial $h$. If $k'\in F[q_-,v,z]$, then $R_kR_{k'}=R_{k+k'}$ for
every polynomial $k$.
\item If $h\in F[q_+,u,z]$ and $k\in F[q_-,v,z]$, then $T_h$ and $R_k$
are automorphisms, with $T_h^{-1}=T_{-h}$ and $R_k^{-1}=R_{-k}$.
\end{enumerate}
\end{lemma}

\begin{proof}
(1) Substituting the images of the generators gives
\[
 T_h(q_+)=(x+2hv+h^2z)z-(v+hz)^2=xz-v^2,\qquad
 R_k(q_-)=(x-2ku-k^2z)z+(u+kz)^2=xz+u^2,
\]
because the terms $2hvz$ and $h^2z^2$, respectively $2kuz$ and $k^2z^2$,
cancel identically. No assumption on $h$ or $k$ is needed. The other
assertions are read off from the tuples.

(2) By (1), $T_h$ fixes every element of $F[q_+,u,z]$; in particular
$T_h(h')=h'$. Hence
\[
 T_hT_{h'}(x)=T_h(x)+2h'\,T_h(v)+h'^2z
 =x+2hv+h^2z+2h'(v+hz)+h'^2z=x+2(h+h')v+(h+h')^2z,
\]
$T_hT_{h'}(v)=T_h(v)+h'z=v+(h+h')z$, and $u$ and $z$ are fixed. The
proof for $R_k$ is the same.

(3) Since $T_0$ is the identity, (2) applied to the pairs $(h,-h)$ and
$(-h,h)$, whose right-hand parameters lie in $F[q_+,u,z]$, gives
$T_hT_{-h}=T_{-h}T_h=1$. The same argument applies to $R_k$.
\end{proof}

The invariance hypothesis in (2) concerns only the right-hand parameter
$h'$, and it cannot be omitted. For example,
$T_uT_x(v)-T_{u+x}(v)=(2uv+u^2z)z\neq0$, because $T_u$ moves the
parameter $x$. Since $\Delta=q_++u^2\in F[q_+,u,z]$ and
$\Delta=q_--v^2\in F[q_-,v,z]$, Lemma~\ref{lem:flows} applies to
$P=T_{\Delta/2}$ and $Q=R_{\Delta/2}$; in particular, $P$ and $Q$ are
automorphisms.

\section{Using a fixed coordinate in a commutator}

Let
\[
 A_+=E_u(q_+/2),\qquad A_-=E_v(q_-/2).
\]
These are elementary. Since $A_+$ and $T_u$ preserve $q_+$,
conjugation replaces the flow parameter $u$ by $u+q_+/2$:
\[
 A_+T_uA_+^{-1}=T_{u+q_+/2}.
\]
The conjugate fixes $u$, because
\[
 A_+T_uA_+^{-1}(u)=A_+T_u(u-q_+/2)=u.
\]
Multiplying by $T_u^{-1}=T_{-u}$ and using Lemma~\ref{lem:flows}(2) leaves
\[
 [A_+,T_u]=T_{q_+/2}.
\]
Thus the auxiliary coordinate is temporarily shifted by the invariant
$q_+/2$, and the commutator produces a flow with precisely that
invariant as its parameter. The same construction yields
$[A_-,R_v]=R_{q_-/2}$.

Using $\Delta=q_++u^2=q_--v^2$ and Lemma~\ref{lem:flows}(2), we obtain
\begin{equation}\label{eq:short-word}
 \boxed{\begin{aligned}
 P&=[E_u(q_+/2),T_u]T_{u^2/2},\\
 Q&=[E_v(q_-/2),R_v]R_{-v^2/2},\\
 \widetilde f&=PQ.
 \end{aligned}}
\end{equation}
For parameters independent of the two moving coordinates, the flows
are products of two elementary maps:
\begin{align*}
 T_h&=E_v(hz)E_x(2hv-h^2z), &&h\in F[u,z],\\
 R_k&=E_u(kz)E_x(-2ku+k^2z), &&k\in F[v,z].
\end{align*}
For example, the first formula sends $x$ to
$x+2h(v+hz)-h^2z=x+2hv+h^2z$. Each block in
\eqref{eq:short-word} therefore expands into
$1+2+1+2+2=8$ elementary factors. The two blocks give sixteen factors
in $(x,v,z,u)$, together with the initial and final linear coordinate
changes when written in the original variables.

The cancellation in the product is visible from
\begin{align*}
 P&=\left(x+\Delta v+\frac{\Delta^2z}{4},
 v+\frac{\Delta z}{2},z,u\right),\\
 Q&=\left(x-\Delta u-\frac{\Delta^2z}{4},
 v,z,u+\frac{\Delta z}{2}\right).
\end{align*}
Since $P$ fixes $\Delta,u,z$, the quadratic correction terms in
$PQ(x)$ cancel, giving $x+\Delta(v-u)$. Both $v$ and $u$ increase
by $\Delta z/2$. Hence $y=v+u$ increases by $\delta z$, whereas
$t=v-u$ is fixed, exactly as required.

The signed permutation $S=(-x,u,z,v)$ gives another group-theoretic
description. It satisfies $S^2=1$, $S(\Delta)=-\Delta$, and
$Q=SP^{-1}S$. Therefore
\[
 \widetilde f=PSP^{-1}S=[P,S].
\]

\section{The Lie algebra and stable tameness}

Assume characteristic zero for the exponential interpretation.
A locally nilpotent derivation $D$ has an exponential
$\exp(D)=\sum_{j\ge0}D^j/j!$, a finite sum on each polynomial.
If $a\in\ker D$, then $aD$ is locally nilpotent as well.

In the original coordinates,
\[
 D_0=t\frac{\partial}{\partial x}+z\frac{\partial}{\partial y},
 \qquad D_0(\delta)=0,\qquad f=\exp(\delta D_0).
\]
The two flows arise from
\[
 D_+=2v\frac{\partial}{\partial x}+z\frac{\partial}{\partial v},
 \qquad
 D_-=-2u\frac{\partial}{\partial x}+z\frac{\partial}{\partial u}.
\]
They are triangular and locally nilpotent, they commute, and they
annihilate $\Delta$.
Under the change to coordinates $(x,v,z,u)$, the derivation $D_0$
becomes
\[
 \widetilde D_0
 =(v-u)\frac{\partial}{\partial x}
 +\frac z2\left(
 \frac{\partial}{\partial v}
 +\frac{\partial}{\partial u}
 \right)
 =\frac{D_++D_-}{2}.
\]
Since $D_+$ and $D_-$ commute and both annihilate $\Delta$, we obtain
\[
 \widetilde f
 =\exp\!\left(\frac{\Delta}{2}D_+\right)
  \exp\!\left(\frac{\Delta}{2}D_-\right)=PQ.
\]
This is the Lie-algebra form of the left-right matrix splitting.

For the commutator defining $P$, put
\[
 X=\frac{q_+}{2}\frac{\partial}{\partial u},\qquad
 Y=uD_+,\qquad Z=\frac{q_+}{2}D_+.
\]
Because $D_+(q_+)=D_+(u)=0$ and $q_+$ is independent of $u$,
\[
 [X,Y]=Z,\qquad [X,Z]=[Y,Z]=0.
\]
These derivations span a three-dimensional Heisenberg Lie algebra.
Consequently
\[
 \exp(X)\exp(Y)\exp(-X)\exp(-Y)=\exp(Z).
\]
Here $\exp(X)=A_+$, $\exp(Y)=T_u$, and
$\exp(Z)=T_{q_+/2}$. The elementary commutator is thus an exact
exponential identity in a two-step nilpotent Lie algebra.

More generally, if $D$ is a triangular locally nilpotent derivation
of $B=F[X_1,\ldots,X_n]$, $a\in\ker D$, and $s$ is a new
variable with $D(s)=0$, the elementary map $A:s\mapsto s+a$
fixing $B$ satisfies
\[
 [A,\exp(sD)]=\exp(aD).
\]
The identity follows by conjugating the parameter $s$ to $s+a$.
Since $A$ is elementary and $\exp(sD)$ is triangular, the extension
of $\exp(aD)$ to $B[s]$ fixing $s$ is tame. Thus $\exp(aD)$ is
stably tame after adjoining one variable.

\smallskip
This is the mechanism behind the stable tameness of the
Nagata automorphism (see the Proposition and the Corollary  in \cite{Smith}); 
an explicit six-factor realization is given in
\cite{Spodzieja}, following the stable-tameness results discussed there.

\smallskip

For the four-variable Anick map, the splitting supplies the auxiliary
coordinates internally: $u$ for the $D_+$-flow and $v$ for the
$D_-$-flow. No fifth variable is introduced. The intermediate
elementary maps may move $\Delta$ and $t$, although the final map
preserves them. All explicit polynomial identities in Sections~2 and~3 have coefficients in $\mathbb Z[1/2]$, so their validity
does not require the exponential notation or characteristic zero.

\section{The sixteen elementary factors}\label{sec:factors}

We work in the coordinates $(x,v,z,u)$, with
\[
 q_+=xz-v^2,\qquad q_-=xz+u^2,\qquad
 \Delta=xz-v^2+u^2.
\]
In the table, $\varepsilon_i=E_w(p)$ means that $w$ is sent to
$w+p$ and every other coordinate is fixed. The factors are ordered so that
\[
 P=\varepsilon_1\cdots\varepsilon_8,\qquad
 Q=\varepsilon_9\cdots\varepsilon_{16},\qquad
 \widetilde f=\varepsilon_1\cdots\varepsilon_{16}.
\]
As throughout the paper, $(AB)(p)=A(B(p))$; thus the rightmost
factor acts first on a polynomial.

\begin{center}
\renewcommand{\arraystretch}{1.3}
\begin{tabular}{@{}rcl@{}}
\toprule
$i$ & Coordinate $w$ & Added polynomial $p$ in $\varepsilon_i=E_w(p)$ \\
\midrule
1  & $u$ & $\tfrac12(xz-v^2)$ \\
2  & $v$ & $uz$ \\
3  & $x$ & $2uv-u^2z$ \\
4  & $u$ & $-\tfrac12(xz-v^2)$ \\
5  & $v$ & $-uz$ \\
6  & $x$ & $-2uv-u^2z$ \\
7  & $v$ & $\tfrac12u^2z$ \\
8  & $x$ & $u^2v-\tfrac14u^4z$ \\
\midrule
9  & $v$ & $\tfrac12(xz+u^2)$ \\
10 & $u$ & $vz$ \\
11 & $x$ & $-2uv+v^2z$ \\
12 & $v$ & $-\tfrac12(xz+u^2)$ \\
13 & $u$ & $-vz$ \\
14 & $x$ & $2uv+v^2z$ \\
15 & $u$ & $-\tfrac12v^2z$ \\
16 & $x$ & $v^2u+\tfrac14v^4z$ \\
\bottomrule
\end{tabular}
\end{center}

\begin{proof}[Completion of the proof of Theorem~\ref{thm:main}]
Each added polynomial omits the coordinate being changed, so all sixteen
maps are elementary; they all fix $z$. Rows~1--8 expand
$A_+T_uA_+^{-1}T_{-u}T_{u^2/2}$, and rows~9--16 expand
$A_-R_vA_-^{-1}R_{-v}R_{-v^2/2}$. Therefore the identities proved in
Section~3 give
\[
 \varepsilon_1\cdots\varepsilon_{16}
 =\left(x+\Delta(v-u),\ v+\frac{\Delta z}{2},\ z,\
 u+\frac{\Delta z}{2}\right).
\]
Returning to $y=v+u$ and $t=v-u$ gives the original Anick
map. The linear coordinate changes also fix $z$ and are tame over
$F[z]$, so the factorization proves tameness over $F[z]$.
The coefficients of the factors and of the linear coordinate changes
belong to $\ZZ[1/2]$.
\end{proof}

\section{The Freudenburg--Moser-Jauslin involution}\label{sec:fmj}

In this section $F$ is a field of characteristic different from $2$
and $3$. We apply the
factorization above to a second automorphism in four variables. Put
\[
 \psi(a,b,x,y)=
 \bigl(a,b,(1-ab)x+a^3y,-b^3x+(1+ab+a^2b^2)y\bigr)
\]
and
\[
 \tau(a,b,x,y)=
 \bigl(b,a,-b^3x+(1+ab+a^2b^2)y,(1-ab)x+a^3y\bigr).
\]
The involution $\tau$ is one of the generators of the action considered
by Freudenburg and Moser-Jauslin in \cite[Section~2]{FMJ}.

\begin{theorem}\label{thm:fmj}
Let $F$ be a field of characteristic different from $2$ and $3$. The
automorphisms $\psi$ and $\tau$ of $F[a,b,x,y]$ are tame. The
automorphism $\psi$ admits a tame factorization in which every factor
fixes $b$; an explicit factorization into $76$ factors, each elementary
with respect to one of three linear coordinate systems, is given in
Appendix~\ref{app:psi}. The factorization identities are
defined over $\ZZ[1/6]$.
\end{theorem}

We retain the composition convention $(AB)(p)=A(B(p))$ from
Section~1. Thus elementary factors written from left to right give
coordinate updates in that order when evaluated as point maps.
We first prove the theorem in characteristic zero, where exponentials
of locally nilpotent derivations are available. Accordingly, from here
on, and in particular in Lemmas~\ref{lem:fmj-separable}
and~\ref{lem:fmj-changes}, $F$ has characteristic zero; this assumption
is dropped in the completion of the proof of Theorem~\ref{thm:fmj},
which extends the result to characteristic different from $2$ and $3$.
The proof starts by expressing $\psi$ as a change of scale of a tame
map obtained from the inverse Anick automorphism. This change of
scale involves $b^{-1}$, so we first adapt the elementary factorization
by separating each factor's value at $b=0$ from its remainder. The
remaining difficulty reduces to four explicit polynomial flows,
which we factor using the commutator construction of Section~3.

\subsection{Comparison with the Anick automorphism}
With $z,t$ renamed $b,a$, the inverse Anick automorphism fixes $a,b$
and acts on the column $(x,y)^{\mathsf T}$ by $C_2$, whereas $\psi$
acts by $C_3$, where
\[
 C_2=\begin{pmatrix}1-ab&a^2\\-b^2&1+ab\end{pmatrix},\qquad
 C_3=\begin{pmatrix}1-ab&a^3\\-b^3&1+ab+a^2b^2\end{pmatrix}.
\]
Indeed, $C_2$ is the inverse of the matrix of
\[
 (x,y,a,b)\longmapsto
 \bigl(x+a(bx-ay),\ y+b(bx-ay),\ a,b\bigr).
\]
For automorphisms fixing $a,b$ that are linear in $x,y$, the matrix of
a product $AB$ is the matrix of $B$ times the matrix of $A$.
At $b=0$, the matrix $C_2$ is elementary. We cancel this specialization
by the elementary shear $U(-a^2)$, with $U$ as in Section~2:
\[
 U(-a^2)=\bigl(C_2|_{b=0}\bigr)^{-1}
  =\begin{pmatrix}1&-a^2\\0&1\end{pmatrix}.
\]
This choice gives
\[
 C_2U(-a^2)=\begin{pmatrix}1-ab&a^3b\\-b^2&1+ab+a^2b^2\end{pmatrix}
 \equiv I\pmod b.
\]
The extra factor $b$ in the upper-right entry suggests rescaling $x$.
After inverting $b$, direct multiplication gives
\begin{equation}\label{eq:fmj-matrix}
 C_3=\Lambda^{-1}C_2U(-a^2)\Lambda,\qquad \Lambda=\operatorname{diag}(b,1).
\end{equation}

Put $v=(y+a)/2$, $u=(y-a)/2$, so that $a=v-u$ and $y=v+u$; this is the
substitution $z\mapsto b$, $t\mapsto a$ in the coordinates of Section~2.
Let $\varepsilon_1,\ldots,\varepsilon_{16}$ be the factors in
Section~\ref{sec:factors}, with $z$ replaced by $b$. Define
\begin{equation}\label{eq:fmj-word}
 e_1=E_x\bigl(-(v-u)^2(v+u)\bigr),\qquad
 e_j=\varepsilon_{18-j}^{-1}\quad(2\leq j\leq17),\qquad
 \Phi=e_1\cdots e_{17}.
\end{equation}
Thus the seventeen factors are specified by the existing table,
with no new choice of factorization; the first factor
$e_1=E_x(-a^2y)$ has the matrix $U(-a^2)$. Direct substitution gives
\begin{equation}\label{eq:fmj-target}
 \begin{gathered}
 \hat x=(1-ab)x+a^3b\,y,\qquad
 \hat y=-b^2x+(1+ab+a^2b^2)y,\\
 \Phi(x,v,b,u)=\bigl(\hat x,(\hat y+a)/2,b,(\hat y-a)/2\bigr).
 \end{gathered}
\end{equation}
In particular, $\Phi$ is the identity modulo $b$.

Let $L$ fix $v,u,b$ and send $x$ to $bx$; its matrix is $\Lambda$.
Under our composition
convention, \eqref{eq:fmj-matrix} says that $L\Phi L^{-1}$ is
$\psi$ in the coordinates $(x,v,b,u)$. Since $L^{-1}$ involves
$b^{-1}$, this identity initially holds in the localization at $b$.
To obtain a tame factorization over $F$, we must arrange that the
rescaled factors are polynomial and tame also at $b=0$.

\subsection{Removing the specialization}
Although $\Phi$ is the identity modulo $b$, its individual elementary
factors need not be. We first separate the part of each factor that
vanishes at $b=0$; this supplies the divisibility needed for the
change of scale.

For each $e_j$, specialize its added polynomial at $b=0$, while
continuing to fix $b$, and denote the resulting elementary map by
$e_j^0$. Put
\[
 r_j=e_j(e_j^0)^{-1},\qquad
 \Phi_j=e_1\cdots e_j,\qquad
 \Phi^0_j=e_1^0\cdots e_j^0,\qquad
 \rho_j=\Phi^0_{j-1}\,r_j\,(\Phi^0_{j-1})^{-1},
\]
with $\Phi_0=\Phi^0_0=1$; thus $\Phi=\Phi_{17}$.
The added polynomial of $r_j$ is obtained by subtracting its value
at $b=0$, and is therefore divisible by $b$. Moving the specialized
factors to the right conjugates each remainder by the preceding
specialized prefix: since $e_j(e_j^0)^{-1}=r_j$,
\[
 \Phi_j(\Phi^0_j)^{-1}=\bigl(\Phi_{j-1}(\Phi^0_{j-1})^{-1}\bigr)\rho_j .
\]
Specialization at $b=0$ is compatible with composition of maps fixing
$b$, and $\Phi$ is the identity modulo $b$; hence $\Phi^0_{17}=1$, and
induction gives
\begin{equation}\label{eq:fmj-relative}
 \Phi=\rho_1\cdots\rho_{17}.
\end{equation}

In the table of Section~\ref{sec:factors}, $x$ enters the added
polynomials only through $xz$, which becomes $xb$. Hence every
specialized factor $e^0_i$ fixes $b$, maps $F[v,u]$ into itself and
sends $x$ to $x$ plus an element of $F[v,u]$, and so does every
$\Phi^0_{j-1}$. Let $\varphi_j$
be the automorphism of $F[v,u]$ induced by $\Phi^0_{j-1}$, extended to
fix $x$ and $b$; it is a product of elementary maps, hence tame. Put
\[
 h_j=x-(\Phi^0_{j-1})^{-1}(x)\in F[v,u].
\]
Applying $\Phi^0_{j-1}$ to $(\Phi^0_{j-1})^{-1}(x)=x-h_j$ gives
$\Phi^0_{j-1}(x)=x+\varphi_j(h_j)$, and therefore
\[
 \Phi^0_{j-1}=\varphi_j\,E_x(h_j).
\]
For $g\in F[v,u]$ we write $E_x(g/b)$ for the translation
$x\mapsto x+g/b$ of the localization $F[b,b^{-1},x,v,u]$; then
$LE_x(g)L^{-1}=E_x(g/b)$. Since $L$ commutes with $\varphi_j$, the
rescaled factor is
\begin{equation}\label{eq:fmj-dilated}
 L\rho_jL^{-1}
 =\varphi_j\,E_x(h_j/b)\,(Lr_jL^{-1})\,E_x(-h_j/b)\,
  \varphi_j^{-1}.
\end{equation}
Thus it suffices to factor the middle expression over the polynomial
ring. The outer maps $\varphi_j$ and $\varphi_j^{-1}$ are already tame.

\subsection{The polynomial factors after the change of scale}\label{subsec:flows}
We now check the middle expression in \eqref{eq:fmj-dilated}.
The rational translation $E_x(h_j/b)$ is only an intermediate map;
the full conjugate will be polynomial. Inspection of the explicit
Anick word \eqref{eq:fmj-word} gives three cases.

If $r_j=E_x(bp)$ with $p\in F[v,u]$, then $Lr_jL^{-1}=E_x(p)$.
It commutes with $E_x(h_j/b)$, so the middle expression is elementary.
This includes $r_1=1$.

If $r_j=E_w(bp)$ for $w\in\{v,u\}$, with $p\in F[v,u]$ independent
of $w$, then $Lr_jL^{-1}=r_j$. Since
$r_j^{-1}E_x(h_j/b)\,r_j=E_x\bigl(r_j^{-1}(h_j)/b\bigr)$, the middle
expression is
\begin{equation}\label{eq:fmj-two-factors}
 r_j\,E_x\!\left(\frac{r_j^{-1}(h_j)-h_j}{b}\right).
\end{equation}
The quotient is a polynomial because $r_j^{-1}$ is the identity
modulo $b$. It does not contain $x$, so both factors are elementary.

The remaining case is $r_j=E_w(\lambda bx)$, where
$w\in\{v,u\}$ and $\lambda\in F$. Its rescaling is
$Lr_jL^{-1}=E_w(\lambda b^2x)$.
Conjugating this shear by $E_x(h_j/b)$ leads to a polynomial flow
whose invariant is $bx+h_j$. We describe it first for $w=v$.

For $h\in F[v,u]$, set
\[
 \theta_h=bx+h(v,u),\qquad
 \mathcal D_h=b\frac{\partial}{\partial v}
       -h_v\frac{\partial}{\partial x}.
\]
The derivation $\mathcal D_h$ is triangular and locally nilpotent,
fixes $u,b$, and annihilates $\theta_h$. For $\lambda\in F$, the
automorphism $N_{h,\lambda}=\exp(\lambda\theta_h\mathcal D_h)$ has the
formula
\begin{equation}\label{eq:fmj-flow}
 \begin{gathered}
 v'=v+\lambda b\theta_h,\qquad
 x'=x-\frac{h(v',u)-h(v,u)}{b},\qquad
 u'=u,\quad b'=b.
 \end{gathered}
\end{equation}
The numerator is divisible by $b$, so this formula is polynomial.
It preserves $\theta_h$, and its inverse is $N_{h,-\lambda}$.
We use the same notation with $u,v$ interchanged when $u$ is the
moving coordinate.

A direct calculation gives
\[
 E_x(h/b)\,E_v(\lambda b^2x)\,E_x(-h/b)=N_{h,\lambda}:
\]
the left side sends $v$ to $v+\lambda b^2(x+h/b)=v+\lambda b\theta_h$
and $x$ to $x+h/b-h(v+\lambda b\theta_h,u)/b$.
Consequently the middle expression in \eqref{eq:fmj-dilated}
is $N_{h_j,\lambda}$ in the remaining case. There are four such
factors in the word. Reversing the preceding
specialized factors gives the following complete list:
\begin{center}
\renewcommand{\arraystretch}{1.3}
\begin{tabular}{@{}rccl@{}}
\toprule
$j$ & $w$ & $\lambda$ & $h_j=x-(\Phi^0_{j-1})^{-1}(x)$\\
\midrule
6  & $v$ & $\tfrac12$  & $vu^2-v^3-2uv-u^3$\\
9  & $v$ & $-\tfrac12$ & $(v-u^2/2)u^2-(v-u^2/2)^3$\\
14 & $u$ & $\tfrac12$  & $2vu-v^3$\\
17 & $u$ & $-\tfrac12$ & $0$\\
\bottomrule
\end{tabular}
\end{center}
The table reduces the proof to four particular flows. Their
polynomiality follows from \eqref{eq:fmj-flow}; their tameness
requires a further argument. The useful feature of these potentials
is that elementary changes reduce them to sums of polynomials in
separate coordinates.

\subsection{Factoring the four remaining flows}
We return to the commutator construction used for the Anick
automorphism. If $h$ separates as $p(r)+q(s)$, the corresponding
derivation splits into commuting parts. For each part, the other
coordinate is available for the commutator.

\begin{lemma}\label{lem:fmj-separable}
Let $F$ have characteristic zero. Suppose that an invertible linear change of $v,u$ gives coordinates
$r,s$ in which $h=p(r)+q(s)$, where $p,q$ are polynomials in one
variable. Then $N_{h,\lambda}$ is tame for every $\lambda\in F$,
with all factors fixing $b$.
\end{lemma}
\begin{proof}
Write $\partial_r,\partial_s$ for the partial derivatives with respect
to the coordinates $(x,r,s,b)$, and put
\[
 \mathcal D_r=b\partial_r-p'(r)\partial_x,\qquad
 \mathcal D_s=b\partial_s-q'(s)\partial_x.
\]
These derivations commute and both annihilate
$\theta_h=bx+p(r)+q(s)$. If $r=\alpha v+\beta u$ and
$s=\gamma v+\eta u$, then
$\mathcal D_h=\alpha\mathcal D_r+\gamma\mathcal D_s$. Hence
\[
 N_{h,\lambda}=\exp(\lambda\theta_h\mathcal D_h)
 =\exp(\lambda\alpha\theta_h\mathcal D_r)
  \exp(\lambda\gamma\theta_h\mathcal D_s).
\]
When $u$ is the moving coordinate,
$\mathcal D_h=\beta\mathcal D_r+\eta\mathcal D_s$ instead, and the
argument is the same.
It suffices to factor $\exp(\mu\theta_h\mathcal D_r)$ for $\mu\in F$;
the other factor is treated in the same way. Since
$\theta_h=\theta_{p(r)}+q(s)$, where $\theta_{p(r)}=bx+p(r)$,
\[
 \exp(\mu\theta_h\mathcal D_r)
 =\exp(\mu\theta_{p(r)}\mathcal D_r)\exp(\mu q(s)\mathcal D_r),
\]
and the second factor is triangular.

The derivation $\mathcal D_r$ fixes $s$, $\theta_{p(r)}$ and $b$, and
its coefficients do not involve $s$. The same commutator calculation as
in Section~3, with $s$ as the auxiliary coordinate, gives
\[
 \bigl[E_s(\mu\theta_{p(r)}),\,\exp(s\mathcal D_r)\bigr]
 =\exp(\mu\theta_{p(r)}\mathcal D_r).
\]
Here $E_s(\mu\theta_{p(r)})$ is elementary and $\exp(s\mathcal D_r)$
is triangular. More explicitly, for $\kappa\in F[s,b]$ the flow
$\exp(\kappa\mathcal D_r)$ is
\[
 r\longmapsto r+b\kappa,\qquad
 x\longmapsto x-\frac{p(r+b\kappa)-p(r)}b,
\]
with $s,b$ fixed, and has the elementary factorization
\[
 \exp(\kappa\mathcal D_r)=E_r(b\kappa)\,
 E_x\!\left(-\frac{p(r)-p(r-b\kappa)}b\right).
\]
Thus the commutator uses the existing coordinate $s$ and no additional
variable. Every factor fixes $b$.
\end{proof}

\begin{lemma}\label{lem:fmj-changes}
Let $F$ have characteristic zero. When $v$ is the moving coordinate, the following operations on $h$
preserve the assertion that every $N_{h,\lambda}$ has a tame
factorization fixing $b$:
\begin{enumerate}
\item adding a polynomial $k(u)$;
\item translating $u$ by a constant;
\item replacing $h(v,u)$ by $h(v+g(u),u)$, where $g\in F[u]$.
\end{enumerate}
The same statements hold with $u,v$ interchanged.
\end{lemma}
\begin{proof}
For the first operation, $\mathcal D_{h+k}=\mathcal D_h$ and
\[
 N_{h+k,\lambda}=N_{h,\lambda}\exp(\lambda k(u)\mathcal D_h).
\]
The two factors commute, and the second is triangular. Its inverse
is triangular as well, so the implication works in both directions.
Each of the other two operations replaces $h$ by $\chi(h)$, where
$\chi$ is the affine translation $u\mapsto u+c$ or the elementary shear
$v\mapsto v+g(u)$, extended to fix $x$ and $b$. Such a map leaves the
moving derivative $\partial_v$ unchanged; it satisfies
$\chi(\theta_h)=\theta_{\chi(h)}$ and
$\chi\mathcal D_h\chi^{-1}=\mathcal D_{\chi(h)}$, hence
$\chi N_{h,\lambda}\chi^{-1}=N_{\chi(h),\lambda}$. All the conjugating
maps and the triangular factors fix $b$.
\end{proof}

Put $f_0(v,u)=vu^2-v^3$. Lemma~\ref{lem:fmj-changes} allows us
to add a polynomial in the fixed coordinate $u$. We therefore seek
a sum $c(v+\beta u)^3+d(v+\eta u)^3$ agreeing with $f_0$
in all terms involving $v$. Matching coefficients gives
\[
 c+d=-1,\qquad c\beta+d\eta=0,\qquad
 3(c\beta^2+d\eta^2)=1.
\]
Taking $\beta=1$ gives the choice $\eta=1/3$, $c=1/2$,
$d=-3/2$. Thus the linear coordinates $r=v+u$, $s=v+u/3$
(so that $\alpha=\gamma=1$ in the notation of
Lemma~\ref{lem:fmj-separable}), whose determinant is $-2/3$, satisfy
\begin{equation}\label{eq:fmj-separation}
 \begin{aligned}
 f_0(v,u)+\frac49u^3
   &=\frac12r^3-\frac32s^3,\\
 f_0(v,u)-v+\frac49u^3
   &=\frac12(r^3+r)-\frac32(s^3+s).
 \end{aligned}
\end{equation}
Lemmas~\ref{lem:fmj-separable} and~\ref{lem:fmj-changes} show that
$N_{f_0,\lambda}$ and $N_{f_0-v,\lambda}$ are tame. The coefficients
lie in $\ZZ[1/6]$, so no extension of $F$ is required.

For row~6 use
\[
 h_6=f_0(v,u-1)-v-u^3.
\]
Translation of $u$ and removal of the term $-u^3$ reduce this case
to $f_0-v$. The potential in row~9 is obtained from $f_0$ by
the shear $v\mapsto v-u^2/2$, so Lemma~\ref{lem:fmj-changes}
applies again.
For row~14, where $u$ is the moving coordinate, remove the term
$-v^3$ and use
\[
 2uv=\frac{(u+v)^2-(u-v)^2}{2},
\]
that is, $r=u+v$ and $s=u-v$ in Lemma~\ref{lem:fmj-separable}.
Finally, row~17 is the elementary map $E_u(-b^2x/2)$.

\begin{proof}[Completion of the proof of Theorem~\ref{thm:fmj}]
First let $F$ have characteristic zero. Every middle expression in
\eqref{eq:fmj-dilated} is polynomial and tame, with all factors fixing
$b$. The same is true of the outer conjugating maps $\varphi_j$. Thus
every $L\rho_jL^{-1}$ has a tame factorization fixing $b$.

Equations~\eqref{eq:fmj-relative} and~\eqref{eq:fmj-dilated} give
\[
 L\Phi L^{-1}=(L\rho_1L^{-1})\cdots(L\rho_{17}L^{-1}).
\]
Both sides are polynomial maps. Since
$F[x,v,b,u]\hookrightarrow F[b,b^{-1},x,v,u]$ is injective,
their equality in the localization is an equality in the polynomial
ring, including at $b=0$. Returning from $v,u$ to $y,a$ uses only
linear changes fixing $b$, and proves the assertion for $\psi$.

Now let $F$ be any field of characteristic different from $2$ and $3$.
Carry out the construction above over $\mathbb Q$. Each factor it
produces is an explicit elementary map: a factor of some $\varphi_j$;
a middle expression $E_x(p)$ or one of the two factors of
\eqref{eq:fmj-two-factors};
the map $E_u(-b^2x/2)$ of row~17; a conjugating map of
Lemma~\ref{lem:fmj-changes}; a map $E_s(\mu\theta_{p(r)})$; or one of
the two elementary factors of a triangular flow, written as for
$\exp(\kappa\mathcal D_r)$ in the proof of
Lemma~\ref{lem:fmj-separable}. Its coefficients are obtained by ring
operations and exact divisions by $b$ from the coefficients of the
table in Section~\ref{sec:factors}, the parameters
$\lambda=\pm\frac12$ and the constants in \eqref{eq:fmj-separation};
hence they lie in $\ZZ[1/6]$. The linear changes of coordinates also
have coefficients in $\ZZ[1/6]$ and determinants that are units of
$\ZZ[1/6]$. The composite of all factors equals $\psi$ over $\mathbb Q$.
Composition of polynomial maps uses only ring operations on the
coefficients, so the same identity holds over $\ZZ[1/6]$, and applying
the homomorphism $\ZZ[1/6]\to F$ to all coefficients gives the identity
over $F$. Each factor remains elementary with respect to its coordinate
system, and each linear change remains invertible. Hence $\psi$ is tame
over $F$, with all factors fixing $b$. The resulting factorization is
listed in Appendix~\ref{app:psi}, where the coefficients can be read
off directly.

Let $\iota(a,b,x,y)=(b,a,y,x)$. In our convention
$\tau=\psi\iota$, so $\tau$ is tame. Furthermore,
\[
 \iota\psi\iota=\psi^{-1},\qquad
 C_3^{-1}=\begin{pmatrix}1+ab+a^2b^2&-a^3\\b^3&1-ab\end{pmatrix}.
\]
Hence $\tau^2=1$, as asserted.
\end{proof}

\section{The \texorpdfstring{$S_3$}{S3}-action}\label{sec:S3}

Take $F=\mathbb C$, let $\omega$ be a primitive cube root of unity,
and put
\[
 \sigma(a,b,x,y)=(\omega a,\omega^{-1}b,x,y).
\]
The associated point maps of $\sigma$ and $\tau$ are exactly the
generators of the action in \cite[Section~2, p.~135]{FMJ}.

\begin{corollary}\label{cor:S3}
The subgroup $\langle\sigma,\tau\rangle$ of
$\TAut(\mathbb C[a,b,x,y])$ is isomorphic to $S_3$ and is not
conjugate to a subgroup of linear automorphisms by any polynomial
automorphism.
\end{corollary}
\begin{proof}
The coefficients $ab,a^3,b^3$ of $\psi$ are invariant under $\sigma$,
so $\psi\sigma=\sigma\psi$. Also
$\iota\sigma\iota=\sigma^{-1}$ and
$\iota\psi\iota=\psi^{-1}$. It follows that
\[
 \sigma^3=\tau^2=1,\qquad
 \tau\sigma\tau=\sigma^{-1}.
\]
These relations show that the group has at most six elements.
Their restrictions $\sigma^j$ and $\sigma^j\tau$, $0\leq j\leq2$,
to $\mathbb C[a,b]$ are distinct, so the group is isomorphic to $S_3$.
The map $\sigma$ is linear, and $\tau$ is tame by
Theorem~\ref{thm:fmj}; hence the group lies in the tame automorphism
group. The algebraic nonlinearizability of this action is
\cite[Theorem~1]{FMJ}.
\end{proof}

The nonlinearizability concerns the whole $S_3$-action.
The cyclic subgroup generated by $\tau$ is linearizable, as noted
in \cite[Section~2]{FMJ}.

\section{Use of AI}\label{sec:AI}

AI tools assisted with developing the factorization,
checking the algebraic identities, and implementing exact symbolic
computations to verify the composition of the sixteen elementary
automorphisms in Section~\ref{sec:factors} and of the $76$ elementary
automorphisms in Appendix~\ref{app:psi}. The authors reviewed all the arguments and computations
and take responsibility for the results presented here.

\appendix
\section{An explicit factorization of \texorpdfstring{$\psi$}{psi}}\label{app:psi}

As an automorphism of $F[a,b,x,y]=F[x,v,b,u]$, where $v=(y+a)/2$ and
$u=(y-a)/2$, the map $\psi$ of Section~\ref{sec:fmj} is determined by
\[
 \psi(x)=(1-ab)x+a^3y,\quad
 \psi(v)=\tfrac12\bigl(\psi(y)+a\bigr),\quad
 \psi(u)=\tfrac12\bigl(\psi(y)-a\bigr),\quad
 \psi(b)=b,
\]
with $a=v-u$, $y=v+u$ and $\psi(y)=-b^3x+(1+ab+a^2b^2)y$.
Besides $\mathcal C_0=(x,v,b,u)$ we use the systems of generators
\[
 \mathcal C_1=(x,r,s,b),\ \ r=v+u,\ s=v+\tfrac13u,
 \qquad\text{and}\qquad
 \mathcal C_2=(x,r,s,b),\ \ r=v+u,\ s=u-v,
\]
so that $v=\frac{3s-r}{2}$, $u=\frac{3(r-s)}{2}$ for $\mathcal C_1$ and
$v=\frac{r-s}{2}$, $u=\frac{r+s}{2}$ for $\mathcal C_2$. A row
$(i,\mathcal C_k,w,p)$ of the table below denotes the automorphism
$\zeta_i$ of $F[x,v,b,u]$ that sends the generator $w$ of
$\mathcal C_k$ to $w+p$ and fixes the other three generators of
$\mathcal C_k$; the polynomial $p$ is written in $\mathcal C_k$ and does
not contain $w$. Thus $\zeta_i$ is elementary with respect to
$\mathcal C_k$, and it fixes $b$. Then
\[
 \psi=\zeta_1\zeta_2\cdots\zeta_{76},
\]
with the convention $(AB)(p)=A(B(p))$ of Section~1; as point maps, the
coordinate updates are performed in the listed order, from row~$1$ to
row~$76$.

The systems $\mathcal C_0$, $\mathcal C_1$ and $\mathcal C_2$ are
obtained from $(a,b,x,y)$ by linear changes of coordinates that fix $b$,
have coefficients in $\ZZ[1/6]$ and have determinants that are units of
$\ZZ[1/6]$. All coefficients in the table also lie in $\ZZ[1/6]$. Hence
the identity holds over every field of characteristic different from
$2$ and $3$, and it exhibits $\psi$ as a product of elementary and
linear automorphisms, all fixing $b$. The rows in $\mathcal C_1$ form
the runs $7$--$22$ and $33$--$48$, and the rows in $\mathcal C_2$ form
the run $55$--$70$. These runs are the words produced by
Lemma~\ref{lem:fmj-separable}, with two commutators and two triangular
flows each, for the flows of rows~$6$, $9$ and $14$
of the table in Section~\ref{subsec:flows}; the remaining rows come from
the outer maps $\varphi_j$, the other middle expressions and
Lemma~\ref{lem:fmj-changes}.

The list was obtained by carrying out the proof of
Theorem~\ref{thm:fmj} literally, which gives $122$ elementary factors,
and merging adjacent factors that move the same generator of the same
system. The identity $\psi=\zeta_1\cdots\zeta_{76}$ was checked by
exact symbolic composition.

{\small
\renewcommand{\arraystretch}{1.25}
\begin{longtable}{@{}rcc>{\raggedright\arraybackslash}p{0.74\textwidth}@{}}
\toprule
$i$ & System & $w$ & Added polynomial $p$\\
\midrule
\endfirsthead
\toprule
$i$ & System & $w$ & Added polynomial $p$\\
\midrule
\endhead
\midrule
\multicolumn{4}{r@{}}{\footnotesize continued on the next page}\\
\endfoot
\bottomrule
\endlastfoot
1 & $\mathcal C_0$ & $x$ & $-\tfrac{1}{4}v^{4}$\\
2 & $\mathcal C_0$ & $u$ & $\tfrac{1}{2}v^{2}b$\\
3 & $\mathcal C_0$ & $x$ & $\tfrac{1}{8}v^{6}b^{2}+\tfrac{1}{4}v^{5}b-\tfrac{3}{4}v^{4}ub-v^{3}u+\tfrac{3}{2}v^{2}u^{2}-v^{2}$\\
4 & $\mathcal C_0$ & $u$ & $vb$\\
5 & $\mathcal C_0$ & $x$ & $v^{3}b^{2}+v^{3}b-3v^{2}ub-2v^{2}u+3vu^{2}+2v^{2}$\\
6 & $\mathcal C_0$ & $u$ & $-1$\\
\midrule
7 & $\mathcal C_1$ & $s$ & $\tfrac{1}{4}r^{3}+\tfrac{1}{2}xb+\tfrac{1}{4}r$\\
8 & $\mathcal C_1$ & $r$ & $sb$\\
9 & $\mathcal C_1$ & $x$ & $-\tfrac{1}{2}s^{3}b^{2}+\tfrac{3}{2}rs^{2}b-\tfrac{3}{2}r^{2}s-\tfrac{1}{2}s$\\
10 & $\mathcal C_1$ & $s$ & $-\tfrac{1}{4}r^{3}-\tfrac{1}{2}xb-\tfrac{1}{4}r$\\
11 & $\mathcal C_1$ & $r$ & $-sb$\\
12 & $\mathcal C_1$ & $x$ & $\tfrac{1}{2}s^{3}b^{2}+\tfrac{3}{2}rs^{2}b+\tfrac{3}{2}r^{2}s+\tfrac{1}{2}s$\\
13 & $\mathcal C_1$ & $r$ & $-\tfrac{3}{4}s^{3}b-\tfrac{3}{4}sb$\\
14 & $\mathcal C_1$ & $x$ & $\tfrac{27}{128}s^{9}b^{2}+\tfrac{81}{128}s^{7}b^{2}+\tfrac{27}{32}rs^{6}b+\tfrac{81}{128}s^{5}b^{2}+\tfrac{27}{16}rs^{4}b+\tfrac{9}{8}r^{2}s^{3}+\tfrac{27}{128}s^{3}b^{2}+\tfrac{27}{32}rs^{2}b+\tfrac{9}{8}r^{2}s+\tfrac{3}{8}s^{3}+\tfrac{3}{8}s$\\
15 & $\mathcal C_1$ & $r$ & $-\tfrac{3}{4}s^{3}+\tfrac{1}{2}xb-\tfrac{3}{4}s$\\
16 & $\mathcal C_1$ & $s$ & $rb$\\
17 & $\mathcal C_1$ & $x$ & $\tfrac{3}{2}r^{3}b^{2}-\tfrac{9}{2}r^{2}sb+\tfrac{9}{2}rs^{2}+\tfrac{3}{2}r$\\
18 & $\mathcal C_1$ & $r$ & $\tfrac{3}{4}s^{3}-\tfrac{1}{2}xb+\tfrac{3}{4}s$\\
19 & $\mathcal C_1$ & $s$ & $-rb$\\
20 & $\mathcal C_1$ & $x$ & $-\tfrac{3}{2}r^{3}b^{2}-\tfrac{9}{2}r^{2}sb-\tfrac{9}{2}rs^{2}-\tfrac{3}{2}r$\\
21 & $\mathcal C_1$ & $s$ & $\tfrac{1}{4}r^{3}b+\tfrac{1}{4}rb$\\
22 & $\mathcal C_1$ & $x$ & $\tfrac{3}{128}r^{9}b^{2}+\tfrac{9}{128}r^{7}b^{2}-\tfrac{9}{32}r^{6}sb+\tfrac{9}{128}r^{5}b^{2}-\tfrac{9}{16}r^{4}sb+\tfrac{9}{8}r^{3}s^{2}+\tfrac{3}{128}r^{3}b^{2}-\tfrac{9}{32}r^{2}sb+\tfrac{3}{8}r^{3}+\tfrac{9}{8}rs^{2}+\tfrac{3}{8}r$\\
\midrule
23 & $\mathcal C_0$ & $v$ & $-\tfrac{2}{9}u^{3}b$\\
24 & $\mathcal C_0$ & $x$ & $-\tfrac{8}{729}u^{9}b^{2}-\tfrac{4}{27}vu^{6}b-\tfrac{2}{3}v^{2}u^{3}+\tfrac{2}{9}u^{5}-\tfrac{2}{9}u^{3}$\\
25 & $\mathcal C_0$ & $u$ & $1$\\
26 & $\mathcal C_0$ & $v$ & $-\tfrac{1}{2}u^{3}b$\\
27 & $\mathcal C_0$ & $x$ & $-\tfrac{1}{8}u^{9}b^{2}-\tfrac{3}{4}vu^{6}b-\tfrac{3}{2}v^{2}u^{3}+\tfrac{1}{2}u^{5}-u^{4}$\\
28 & $\mathcal C_0$ & $v$ & $\tfrac{1}{2}u^{2}$\\
29 & $\mathcal C_0$ & $x$ & $-v^{2}$\\
30 & $\mathcal C_0$ & $u$ & $-vb$\\
31 & $\mathcal C_0$ & $x$ & $\tfrac{1}{8}v^{6}b^{5}+\tfrac{3}{4}v^{5}ub^{4}+\tfrac{15}{8}v^{4}u^{2}b^{3}-\tfrac{3}{4}v^{5}b^{3}+\tfrac{5}{2}v^{3}u^{3}b^{2}-3v^{4}ub^{2}-\tfrac{1}{2}v^{4}b^{3}+\tfrac{15}{8}v^{2}u^{4}b-\tfrac{9}{2}v^{3}u^{2}b-2v^{3}ub^{2}+\tfrac{3}{4}vu^{5}+\tfrac{3}{2}v^{4}b-3v^{2}u^{3}-3v^{2}u^{2}b+3v^{3}u+v^{3}b-2vu^{3}+2v^{2}u$\\
32 & $\mathcal C_0$ & $v$ & $-\tfrac{1}{2}u^{2}$\\
\midrule
33 & $\mathcal C_1$ & $s$ & $-\tfrac{1}{4}r^{3}-\tfrac{1}{2}xb$\\
34 & $\mathcal C_1$ & $r$ & $sb$\\
35 & $\mathcal C_1$ & $x$ & $-\tfrac{1}{2}s^{3}b^{2}+\tfrac{3}{2}rs^{2}b-\tfrac{3}{2}r^{2}s$\\
36 & $\mathcal C_1$ & $s$ & $\tfrac{1}{4}r^{3}+\tfrac{1}{2}xb$\\
37 & $\mathcal C_1$ & $r$ & $-sb$\\
38 & $\mathcal C_1$ & $x$ & $\tfrac{1}{2}s^{3}b^{2}+\tfrac{3}{2}rs^{2}b+\tfrac{3}{2}r^{2}s$\\
39 & $\mathcal C_1$ & $r$ & $\tfrac{3}{4}s^{3}b$\\
40 & $\mathcal C_1$ & $x$ & $-\tfrac{27}{128}s^{9}b^{2}+\tfrac{27}{32}rs^{6}b-\tfrac{9}{8}r^{2}s^{3}$\\
41 & $\mathcal C_1$ & $r$ & $\tfrac{3}{4}s^{3}-\tfrac{1}{2}xb$\\
42 & $\mathcal C_1$ & $s$ & $rb$\\
43 & $\mathcal C_1$ & $x$ & $\tfrac{3}{2}r^{3}b^{2}-\tfrac{9}{2}r^{2}sb+\tfrac{9}{2}rs^{2}$\\
44 & $\mathcal C_1$ & $r$ & $-\tfrac{3}{4}s^{3}+\tfrac{1}{2}xb$\\
45 & $\mathcal C_1$ & $s$ & $-rb$\\
46 & $\mathcal C_1$ & $x$ & $-\tfrac{3}{2}r^{3}b^{2}-\tfrac{9}{2}r^{2}sb-\tfrac{9}{2}rs^{2}$\\
47 & $\mathcal C_1$ & $s$ & $-\tfrac{1}{4}r^{3}b$\\
48 & $\mathcal C_1$ & $x$ & $-\tfrac{3}{128}r^{9}b^{2}-\tfrac{9}{32}r^{6}sb-\tfrac{9}{8}r^{3}s^{2}$\\
\midrule
49 & $\mathcal C_0$ & $v$ & $\tfrac{2}{9}u^{3}b$\\
50 & $\mathcal C_0$ & $x$ & $\tfrac{8}{729}u^{9}b^{2}-\tfrac{4}{27}vu^{6}b+\tfrac{2}{3}v^{2}u^{3}-\tfrac{2}{9}u^{5}+\tfrac{1}{4}u^{4}$\\
51 & $\mathcal C_0$ & $v$ & $-\tfrac{1}{2}u^{2}b$\\
52 & $\mathcal C_0$ & $x$ & $-\tfrac{1}{8}u^{6}b^{2}-\tfrac{3}{4}vu^{4}b-\tfrac{3}{2}v^{2}u^{2}+u^{2}$\\
53 & $\mathcal C_0$ & $v$ & $ub$\\
54 & $\mathcal C_0$ & $x$ & $u^{3}b^{2}-3vu^{2}b+3v^{2}u-2u^{2}$\\
\midrule
55 & $\mathcal C_2$ & $s$ & $\tfrac{1}{2}xb+\tfrac{1}{4}r^{2}$\\
56 & $\mathcal C_2$ & $r$ & $sb$\\
57 & $\mathcal C_2$ & $x$ & $\tfrac{1}{2}s^{2}b-rs$\\
58 & $\mathcal C_2$ & $s$ & $-\tfrac{1}{2}xb-\tfrac{1}{4}r^{2}$\\
59 & $\mathcal C_2$ & $r$ & $-sb$\\
60 & $\mathcal C_2$ & $x$ & $\tfrac{1}{2}s^{2}b+rs$\\
61 & $\mathcal C_2$ & $r$ & $-\tfrac{1}{4}s^{2}b$\\
62 & $\mathcal C_2$ & $x$ & $\tfrac{1}{32}s^{4}b+\tfrac{1}{4}rs^{2}$\\
63 & $\mathcal C_2$ & $r$ & $\tfrac{1}{2}xb-\tfrac{1}{4}s^{2}$\\
64 & $\mathcal C_2$ & $s$ & $rb$\\
65 & $\mathcal C_2$ & $x$ & $-\tfrac{1}{2}r^{2}b+rs$\\
66 & $\mathcal C_2$ & $r$ & $-\tfrac{1}{2}xb+\tfrac{1}{4}s^{2}$\\
67 & $\mathcal C_2$ & $s$ & $-rb$\\
68 & $\mathcal C_2$ & $x$ & $-\tfrac{1}{2}r^{2}b-rs$\\
69 & $\mathcal C_2$ & $s$ & $\tfrac{1}{4}r^{2}b$\\
70 & $\mathcal C_2$ & $x$ & $-\tfrac{1}{32}r^{4}b+\tfrac{1}{4}r^{2}s$\\
\midrule
71 & $\mathcal C_0$ & $u$ & $-\tfrac{1}{2}v^{3}b$\\
72 & $\mathcal C_0$ & $x$ & $v^{4}$\\
73 & $\mathcal C_0$ & $u$ & $-\tfrac{1}{2}v^{2}$\\
74 & $\mathcal C_0$ & $x$ & $u^{2}$\\
75 & $\mathcal C_0$ & $v$ & $-ub$\\
76 & $\mathcal C_0$ & $u$ & $-\tfrac{1}{2}xb^{2}+\tfrac{1}{2}v^{2}$\\
\end{longtable}
}

\end{document}